\documentclass[12pt,a4paper]{amsart}
\usepackage{amssymb}
\usepackage{amsmath}

\usepackage[T1]{fontenc}
\usepackage[utf8]{inputenc}
\usepackage[german,english]{babel}

\usepackage[babel,german=guillemets]{csquotes}
\usepackage[backend=bibtex8,style=authoryear]{biblatex}
\bibliography{piaget.bib}
\usepackage{lmodern}

\usepackage[bottom]{footmisc}

\begin{document}
\pagestyle{plain}

\author{Hans Joachim Burscheid $\cdot$ Horst Struve}
\address{Universität zu Köln, Seminar für Mathematik und ihre Didaktik \\ Gronewaldstr.2, 50931 Köln \\ Tel.: 0049-221-470-4750, Fax.: 0049-221-470-4985 \\ e-mail: h.burscheid@uni-koeln.de, h.struve@uni-koeln.de}

\title{Die Bedeutung empirischer Theorien für die Entwicklung von Gruppierungen}

\selectlanguage{german}
\begin{abstract}
Bei der Rezeption der Denkpsychologie Piagets durch die Mathematikdidaktik der Bundesrepublik Deutschland Ende der 1960er Jahre wurde der Begriff der Gruppierung ähnlich einer mathematischen Struktur aufgefaßt. Kann die Berücksichtigung empirischer Theorien zu einem adäquateren Verständnis beitragen? 
\end{abstract}

\maketitle
\vspace{-1cm}

\selectlanguage{english}
\begin{abstract}
Reviewing Piaget's psychology of reasoning by the mathematical educators of the Federal Republic of Germany late in the 1960$^{\text{th}}$ the concept of grouping has been understood similar to a mathematical structure. Possibly to consider empirical theories can contribute to an understanding more adequate. 
\end{abstract}

\vspace{5mm}
\thispagestyle{empty}

\noindent Als Ende der 1960er Jahre die Rezeption der Piagetschen Denkpsychologie durch die Mathematikdidaktik der Bundesrepublik Deutschland in voller Blüte stand, faszinierte vor allem der Begriff der \emph{Gruppierung}. Die Elemente einer Gruppierung --- die \emph{Operationen} oder \emph{verinnerlichten Handlungen} --- erfüllen folgende Bedingungen: \\

\begin{small}
1. Komposition: x + x$^{\prime}$ = y, y + y$^{\prime}$ = z; etc.

2. Reversibilität: y - x = x$^{\prime}$ oder y - x$^{\prime}$ = x

3. Assoziativität: (x + x$^{\prime}$) + y$^{\prime}$ = x + (x$^{\prime}$ + y$^{\prime}$) = (z).

4. Allgemeine, identische Operation: x - x = o; y - y = o; etc.
 
5. Tautologie oder besondere, identische Operation: x + x = x; y + y = \\ 
\noindent \hspace*{9mm}y; etc. (Piaget [1971, S. 49])  \\
\end{small}

\noindent Etwas anders formuliert: \\

\begin{small}
1. \& 2. x + x$^{\prime}$ = y $\leftrightarrow$ y + (-x) = x$^{\prime}$ $\vee$ y + (-x$^{\prime}$) = x

\noindent \hspace*{13mm}3. (x + x$^{\prime}$) + y$^{\prime}$ = x + (x$^{\prime}$ + y$^{\prime}$)

\noindent \hspace*{13mm}4. x + (-x) = (-x) + x = o

\noindent \hspace*{13mm}5. x + x = x \\
\end{small}

In dieser Formulierung sieht man deutlicher, daß es nur eine Verknüpfung (+) für die Operationen gibt und daß es bzgl. der allgemeinen, identischen Operation (o) zu gewissen Operationen x je eine inverse Operation -x geben kann.   

Der Piagetsche Text legte es nahe, den Begriff der Gruppierung so zu verstehen, als besäße eine Gruppierung Modelle, in denen jede dieser Bedingungen erfüllt sei. Heinz Griesel wies als erster darauf hin, daß eine Gruppierung bei dieser Lesart nur einelementig sein könne, also nur ein triviales Modell habe [1970, S. 128].

Wenn man heute --- rund 40 Jahre später --- diese Diskussion noch einmal nachliest, stellt man erstaunt fest, wie selbstverständlich man annahm, Piagets Formulierungen als mathematisch tragfähig ansehen zu können. Eine Sichtweise, die allerdings nicht auf deutsche Mathematikdidaktiker beschränkt war.

\begin{quote}
\begin{small}
\enquote{A grouping incorporates properties from two well--known mathematical structures, the \emph{group} and the \emph{lattice}. A Piagetian grouping is thus a kind of hybrid group--lattice structure.} (Flavell [1985, S. 92])
\end{small}
\end{quote}

Mathematisch tragfähig erschienen Piagets Formulierungen in dem Sinne, daß die im folgenden genannten Kollegen offenbar davon ausgingen, daß die von ihm angegebenen Bedingungen in der Tat wie ein mathematisches Axiomensystem gelesen werden könnten. Dazu dürfte die Bezeichnung \enquote{Gruppierung} --- in Analogie zur (mathematischen) \enquote{Gruppe} --- und die \emph{formale} Darstellung der Bedingungen  entscheidend beigetragen haben, die aber vermutlich in erster Linie Piagets Neigung zur Mathematik zu danken waren. Denn die Bedingungen wie ein Axiomensystem zu lesen wirft doch einige wichtige Fragen auf.

\sloppy Erich Wittmann [1973] gab als erster das Axiomensystem einer Struktur an, in der --- unter Berücksichtigung der Grieselschen Feststellung --- die Bedingungen von Piaget hineininterpretiert werden konnten. Um die Bedingung der Tautologie zu erfassen führt er neben + eine weitere Verknüpfung ein, indem er mit Hilfe von + eine Quasiordnung unter den Operationen definiert, an die er verbandstheoretische Forderungen stellt. Hans -- Georg Steiner [1973] hat das von Wittmann angegebene Axiomensystem aufgegriffen und in eine vereinfachte Form gebracht. Elmar Cohors -- Fresenborg [1974] beschrieb schließlich die Wittmannsche Gruppierung als ein Semi -- Thue -- System. Die von diesen Autoren angegebenen Beispiele für Gruppierungen sind zeitbedingt stark an der \enquote{Neuen Mathematik} orientiert und wirken heute etwas artifiziell. Natürlich genügen sie dem Axiomensystem und eine gewisse Komplexität der Beispiele ist zwar einleuchtend, da Gruppierungen als Gleichgewichtsform der Operationen/verinnerlichten Handlungen den "Endzustand" der Stufe des konkret--operationalen Denkens beschreiben, auf der sich die Schüler\footnote{und natürlich auch die Schülerinnen} befanden, für die die Beispiele gedacht waren. Doch erscheint es wenig einleuchtend, daß eine derart komplizierte Gleichgewichtsform, wie sie die Wittmannsche Struktur beschreibt, ohne eine sehr spezifische Förderung erreicht werden kann. Und dies überzeugt nicht. Denn nach der Auffassung Piagets erwerben wir die von ihm formulierten Fähigkeiten im Alltag. Verinnerlicht werden  Handlungen, die real durchgeführt werden, eben weit überwiegend im Alltag und nicht in einer Unterrichtssituation. Es sind in der Regel solche, die auf die WELT ausgeübt werden. 

Die   Operationen/verinnerlichter Handlungen werden von Piaget wie folgt charakterisiert:

\begin{quote}
\begin{small}
\enquote{Die spezifische Natur der Operationen besteht, verglichen mit den empirischen Tätigkeiten, gerade in der Tatsache, daß sie niemals in diskontinuierlichem Zustand existieren. Es ist nur eine gänzlich unerlaubte Abstraktion, wenn man von \flqq einer\frqq $\:$Operation spricht; eine vereinzelte Operation kann nicht Operation sein, denn die eigentümlichste Eigenschaft der Operationen liegt gerade darin, daß sie zu Systemen vereinigt sind.} [1971, S. 41]  
\end{small}
\end{quote}

Es stellt sich die Frage nach der Organisationsform der real durchgeführten Handlungen, die den Operationen zugrunde liegen. Denn man kann zwar eine real durchgeführte Handlung als einzelne betrachten, aber aus einer einzelnen Handlung läßt sich kein \enquote{operatives Gesamtsystem} (Piaget) ableiten und aus einer Vielzahl von planlos durchgeführten Handlungen ebenfalls nicht.

Das Verhalten von Kindern und damit auch ihr Handeln ist stark regelhaft. Sie verhalten sich mitunter, als verfügten sie über gewisse Theorien, die ihr Handeln steuern (Gopnik und Meltzoff [1997]). Eine solche Theorie könnte z.B. ihren Umgang mit einem mit Luft gefüllten Ball oder mit einem Luftballon steuern. Diese Theorien beziehen sich unmittelbar auf die WELT, ihre Begriffe haben eine starke ontologische Bindung. Man spricht gemeinhin von \emph{empirischen Theorien}. Eine etablierte Form zur \emph{Darstellung} empirischer Theorien ist die strukturalistische Metatheorie, die von Wolfgang Stegmüller und seinem Kreis entwickelt wurde [1973, 1986], Balzer [1982]. Ein wesentliches Merkmal dieser Darstellungsform ist, daß sie den Aufbau der Theorie wiedergibt.

Die Darstellung gliedert sich in drei Stufen. Auf der ersten Stufe werden als sog. \emph{partielle Modelle} der Theorie die empirischen Gegebenheiten formuliert, über die die Theorie Aussagen machen kann. Diese werden aus sog. \emph{paradigmatischen Beispielen} abgeleitet. Im vorliegenden Fall betrifft dies die Objekte, auf die reale Handlungen ausgeübt werden sowie diese Handlungen selbst. Die nächste Stufe --- die der \emph{potentiellen Modelle} --- betrifft die Sprache, in der die Theorie formuliert wird. Sie soll möglichst präzise und deshalb formalisierbar sein. Es liegt daher nahe, auf Begriffe der Mathematik, z.B. Relationen oder Funktionen, zurückzugreifen. Sie erlauben es, die partiellen Modelle formal zu beschreiben. Dies ist in dem Sinne von besonderem Interesse als auf dieser Stufe die realen Handlungen der ersten Stufe durch mathematische Strukturen beschrieben, d.h. in Systeme eingebunden werden. Die Systeme der realen Handlungen bilden die Vorlage für die Systeme der verinnerlichten Handlungen. Während auf der ersten Stufe noch alle Begriffe reale Referenzen haben oder aus bekannten (früher erworbenen) Theorien stammen, bedarf die neu zu formulierende Theorie neuer Begriffe, will sie neues Wissen vermitteln. Es sind die bzgl. der in den partiellen Modellen formulierten empirischen Gegebenheiten \emph{theoretischen Begriffe}, solche, die dort keine Referenzen haben, die erst durch die neu zu formulierende Theorie eine Bedeutung erhalten. Um diese wird die Sprache der Theorie erweitert. Auf der dritten Stufe werden schließlich die Axiome formuliert, die die \emph{Modelle} der Theorie definieren. Man kann die Modelle als die formalen Spiegelbilder der Systeme der Operationen ansehen, die auf der kognitiven Ebene das Verfügen über die empirische Theorie ausmachen. Dies will sagen, daß eine empirische Theorie zu erwerben --- psychologisch betrachtet ---  bedeutet, bestimmte Handlungen zu verinnerlichen.

Die Organisationsstruktur der verinnerlichten Handlungen, die man entwickelt, wenn man über eine empirische Theorie zu verfügen lernt, kann einzelne der Piagetschen Bedingungen erfüllen. Alle Bedingungen einer Gruppierung gleichzeitig zu erfüllen ist --- in der Piagetschen Lesart --- wegen der Unvereinbarkeit von allgemeiner und besonderer, identischer Operation ohnehin nicht möglich.

Die Eigenschaften einer Gruppierung sind bei dieser Sichtweise ein Extrakt der Eigenschaften real durchgeführter Handlungen, die der einzelne verinnerlicht, wenn er über unterschiedliche empirische Theorien zu verfügen lernt. Ein so verstandener Gruppierungsbegriff enzieht sich einer mathematischen Modellierung, trifft aber die psychologischen Intentionen, die Piaget verfolgt hat. Die empirischen Theorien, über die Kinder zu verfügen lernen, sind die natürlichen Zugänge zur Entwicklung von Gruppierungen.

Es kann durchaus sein, daß einzelne der Piagetschen Bedingungen auch in der Formulierung der Modelle auftreten, aber dies wäre eine rein sprachliche Übereinstimmung, denn die Objekte, von denen die Modelle sprechen, sind nicht die Operationen Piagets und die verwendeten Begriffe sind nicht notwendig identisch mit den seinen. Dies zeigt nicht nur die Beobachtung Griesels. So betrachtet Piaget z.B. die \emph{Gruppierung der asymmetrischen Relationen} und schreibt:
 
\begin{quote}
\begin{small}
\enquote{Nennen wir a die Relation O < A; b die Relation O < B; c die Relation O < C. Man kann dann die Relation A < B a$^{\shortmid}$ nennen, die Relation B < C b$^{\shortmid}$ etc. Die umgekehrte Operation besteht aus der Subtraktion einer Relation, was der Addition ihrer Konverse äquivalent ist. ..... Auf der Transitivität, die dieser Seriation eigentümlich ist, gründet sich der Schluß A < B; B < C, also A < C.} [1971, S. 51]
\end{small}
\end{quote}
 
Es ergibt sich folgende Konsequenz:

\begin{center}
\begin{small}
a + (-a) = (O < A $\wedge$ A < O) = (O < O)
\end{small} 
\end{center} 

Eine (mathematische) asymmetrische Relation ist irreflexiv, was die Piagetsche offensichtlich nicht ist. Nur aus der sprachlichen Übereinstimmung auch auf eine inhaltliche zwischen Piagets Begriffen und den Begriffen der Mathematik zu schließen, dürfte wenig zweckdienlich sein.  

Wie wir wiederholt betont haben, sind wir der Auffassung, daß die Inhalte der Elementarmathematik von Schülern im Rahmen empirischer Theorien erworben werden (Ein Standpunkt, der z.B. auch von Griesel geteilt wird [2013]). Elementarmathematik umfaßt dabei mindestens die Inhalte, die auf der konkret -- operationalen Stufe oder zuvor erlernt werden, also i.w. alle Inhalte der Klassen 1 bis 10. Ein zentrales Anliegen des Unterrichts sollte es daher sein, diese Inhalte in empirische Theorien zu integrieren. Für die Zahlbegriffsentwicklung haben wir an etlichen Beispielen gezeigt, wie dies zu verstehen ist [2009]. Ein Unterricht, der diese Sichtweise berücksichtigt, trägt ganz selbstverständlich zur Ausbildung von Gruppierungen als Gleichgewichtszuständen des konkret--operationalen Denkens bei. Den Kindern solche wie \enquote{isomorphe Spielhandlungen} (Breidenbach) vorzustellen, weist sicherlich nicht den Weg, auf dem sie diese entwickeln (vgl. Bussmann [1974]). \medskip

\noindent \emph{Bemerkung}: Verzichtet man darauf --- wie hier vorgeschlagen --- die Piagetschen Bedingungen wie ein Axiomensystem zu lesen, so erledigt sich die Unvereinbarkeit von allgemeiner (4.) und identischer Operation (5.), auf die Griesel hingewiesen hatte. Während in Prozessen, die sich auf --- im weitesten Sinne --- quantifizierbare Objekte beziehen, einzelne oder alle der Bedingungen 1. -- 4. erfüllt sein mögen, gilt dies nicht unbedingt für solche, die sich auf Qualitäten beziehen. Betrachtet man z.B. den Prozeß des Rührens eines Teiges, so ist dieser nicht umkehrbar und damit entfällt 4. Unterbricht man das Rühren und rührt erneut, so ruft dies keine Änderung hervor, d.h. diese Tätigkeit ist in der Sprache Piagets tautologisch. Daß Piaget solche Beispiele vor Augen gehabt haben muß, belegt folgendes Zitat:

\begin{quote}
\begin{small}
\enquote{Im Bereich der Zahlen bildet eine zu sich selbst addierte Einheit eine neue Zahl durch Anwendung der Komposition ( ... ). Es findet eine Iteration statt. Im Gegensatz dazu verändert sich ein qualitatives Element nicht durch Wiederholung, sondern ergibt eine \flqq Tautologie \frqq: A + A = A.} [1971, S. 48/49]
\end{small}
\end{quote}

Aus der Formulierung wird deutlich, daß Piaget, was er ein qualitatives Element nennt, als etwas nicht Quantifizierbares auffaßt. Die Elemente einer Gruppierung sind aber Operationen/verinnerlichte Handlungen. Die Vermutung liegt nahe, daß er solche verinnerlichten Handlungen meinte, deren Realisierungen sich auf Qualitäten beziehen.     

\selectlanguage{german}
\nocite{*}
\printbibliography

\end{document}